\documentclass{article}
\usepackage[american]{babel}
\usepackage{amsfonts,amsmath,amssymb,epsf,epsfig}
\usepackage{xypic}
\xyoption{all}
\input{xypic}
\newdir{ >}{{}*!/-8pt/\dir{>}}
\def\dar[#1]{\ar@<2pt>[#1]\ar@<-2pt>[#1]}

\def\R{\mathbb{R}}

\def\C{\mathbb{C}}

\def\Z{\mathbb{Z}}

\def\coker{{\mathrm coker}}
\def\ker{{\rm ker}}
\def\Alt{{\rm Alt}}

\def\coker{\mathrm coker}
\def\id{\rm id}

\def\pr{\noindent $\bf{Proof.}$\quad}     
\def\fin{\hfill$\square$\\}           
\newtheorem{theo}{Theorem}
\newtheorem{defi}{Definition}
\newtheorem{rem}{Remark}
\newtheorem{prop}{Proposition}
\newtheorem{cor}{Corollary}

\newtheorem{lem}{Lemma}
\begin{document}
\title{Obstruction classes of crossed modules of Lie algebroids
and Lie groupoids linked to existence of principal bundles} 

\author{Camille Laurent-Gengoux\\
        Universit\'e de Poitiers\\
\and    Friedrich Wagemann \\
        Universit\'e de Nantes}
   
\maketitle


\begin{abstract}
Let $K$ be a Lie group and $P$ be a $K$-principal bundle on a manifold 
$M$. Suppose given furthermore a central extension 
$$1\to Z\to \hat{K}\to K\to 1$$
of $K$. It is a classical question whether there exists a $\hat{K}$-principal 
bundle $\hat{P}$ on $M$ such that $\hat{P}/Z\cong P$.
In \cite{Neeb}, Neeb defines in this context a crossed module of topological 
Lie algebras whose cohomology class $[\omega_{\rm top\,\,alg}]$ is an 
obstruction to the existence of $\hat{P}$. 
In the present paper, we show that $[\omega_{\rm top\,\,alg}]$ is up to 
torsion a full obstruction for this problem, and we clarify its relation
to crossed modules of Lie algebroids and Lie groupoids, and finally to gerbes. 
\end{abstract}

\section*{Introduction}

It is well known that $3$-cohomology in an algebraic category like the 
categories of (discrete) groups, Lie algebras or associative algebras is 
related to crossed modules of groups, Lie algebras or associative algebras
respectively. A {\it crossed module of groups} is, roughly speaking, a 
homomorphism of groups $\mu:M\to N$ together with an action by automorphisms 
of $N$ on $M$ which is compatible in some sense with $\mu$. Passing to 
kernel and cokernel of $\mu$, one gets a $4$-term exact sequence of groups
$$1\to V=\ker\,\mu\to M\to N\to G=\coker\,\mu\to 1,$$ 
such that $V$ is an abelian group and a $G$-module. The general algebraic 
picture associates to such a crossed module a cohomology class in $H^3(G,V)$
which is the obstruction to come from an extension of $G$ by $M$. 
Related notions which take into account topology exist for Lie groupoids 
\cite{Mack1}, \cite{Andr} and for topological Lie algebras \cite{Neeb}, 
\cite{Wage}. 

In \cite{Neeb}, Karl-Hermann Neeb defines for a crossed module of topological 
Lie algebras which is split as a sequence of topological vector spaces a
cohomology class $[\omega_{\rm top\,\,alg}]$. He shows that 
$[\omega_{\rm top\,\,alg}]$ has a specific meaning in the following context:
let $K$ be some Lie group and $P$ be a $K$-principal bundle on some manifold 
$M$. Suppose given furthermore a central extension 
$$1\to Z\to \hat{K}\to K\to 1$$
of $K$. The question is now whether there exists a $\hat{K}$-principal 
bundle $\hat{P}$ on $M$ such that $\hat{P}/Z\cong P$. Neeb uses the 
ingredients of the problem to associate a crossed module of topological Lie 
algebras to it such that its obstruction class $[\omega_{\rm top\,\,alg}]$
is a $3$-de Rham cohomology class on $M$ which is an obstruction to the 
existence of $\hat{P}$. 

The origin of this paper is the question of Neeb at the end of section IV 
of \cite{Neeb} whether $[\omega_{\rm top\,\,alg}]$ is the full obstruction 
to the existence of $\hat{P}$, which $3$-de Rham forms arise as obstructions,
and on the relation of $[\omega_{\rm top\,\,alg}]$ to gerbes.

We answer his questions in the framework of crossed modules of Lie algebroids 
and groupoids and show that $[\omega_{\rm top\,\,alg}]$ can be identified 
with the obstruction class of a certain crossed module of Lie algebroids
associated to the above problem (theorem $1$), and, up to torsion, even to 
the obstruction class of a certain crossed module of Lie groupoids
associated to the above problem (theorem $2$, main theorem of this paper), 
which is known \cite{Mack1} to be the full obstruction to the existence of 
$\hat{P}$.

In section $5$, we show that it follows from Serre's identification of the 
Brauer group $Br(M)$ of $M$ (cf \cite{Grot2}) that Neeb's obstruction class
is zero for finite dimensional structure group. In section $7$, we deduce from
the observation that gerbes and crossed modules of Lie groupoids are classified
by the same kind of cohomology classes a direct relation between these two 
kinds of objects.

In the appendix on Deligne cohomology, we present the 
necessary material for the proof of theorem $2$. 

{\bf Acknowledgements:} We are indepted to Karl-Hermann Neeb for spotting and 
correcting an error. Both authors thank the Rencontres
Math\'ematiques de Glanon where some of the work presented here was done.

\section{The Atiyah sequence}

In this section, we recall the Atiyah sequence associated to a $K$-principal 
bundle on a fixed base manifold $M$, and we explain the main object of this 
article.  

Let $M$ be a (finite dimensional, connected, second countable, smooth) 
manifold with finite dimensional de Rham cohomology, 
and let $K$ be a (not necessarily finite dimensional) Lie group 
with Lie algebra ${\mathfrak k}$. We usually take infinite dimensional 
Lie groups to be modeled on locally convex spaces. 
Furthermore, let $\pi:P\to M$ be a $K$-principal bundle on $M$. As for finite 
dimensional structure groups, $P$ can be represented by a smooth \v{C}ech 
$1$-cocycle, cf \cite{Wock}. Connections on $P$ can be constructed by 
patching on the finite dimensional manifold $M$. $K$ acts on 
$P$ from the right, and this action induces an action on $TP$. The induced 
map between tangent bundles $T\pi:TP\to TM$ factors to a map 
$\pi_*:TP/K\to TM$. The kernel of $\pi_*$ can be identified (cf \cite{Mack} 
p. 92 or \cite{Neeb} IV) with $(P\times{\mathfrak k})/K$, where $K$ acts on 
the product by the diagonal action, using the adjoint action on the second 
factor. We denote this bundle by ${\rm Ad}(P)$. Therefore, we get the 
{\it Atiyah sequence} of vector bundles
$$0\to {\rm Ad}(P)\to TP/K\to TM \to 0.$$ 

The main question we address in this paper is the following: 
given a central extension 
$$1\to Z\to \hat{K}\to K\to 1$$
of $K$ by $Z$, is there a $\hat{K}$-principal bundle $\hat{P}$ on $M$ such 
that $\hat{P}/Z\cong P$ ? More precisely, one wants to construct computable
obstructions for the existence of such a $\hat{P}$. This is the point of view
expressed in Neeb's paper \cite{Neeb} section VI.   

\section{Crossed modules of topological Lie algebras}

In \cite{Neeb}, Neeb associates to a given principal bundle $P$ on $M$ 
and a central extension 
$$1\to Z\to \hat{K}\to K\to 1$$
of its structure group $K$, a differential $3$-form $\omega_{\rm top\,\,alg}$ 
defining a class in $H^3_{\rm dR}(M,{\mathfrak z}) =H^3_{\rm dR}(M)
\otimes {\mathfrak z} $ which is an obstruction to 
the existence of a principal $\hat{K}$-bundle $\hat{P}$ such that 
$\hat{P}/Z\cong P$. Let us recall its construction:

\begin{defi}
A crossed module of Lie algebras is a homomorphism of
Lie algebras $\mu:\mathfrak{m}\to\mathfrak{n}$ together with an action $\eta$ 
of $\mathfrak{n}$ on $\mathfrak{m}$ by derivations such that

(a) $\mu(\eta(n)\cdot m)\,=\,[n,\mu(m)]$ for all $n\in\mathfrak{n}$
and all $m\in\mathfrak{m}$,

(b) $\eta(\mu(m))\cdot m'\,=\,[m,m']$ for all $m,m'\in\mathfrak{m}$.
\end{defi}

In the framework of {\it topological Lie algebras}, one requires all maps to 
be continuous and topologically split, cf \cite{Neeb}. 

To an Atiyah sequence
$$0\to {\rm Ad}(P)\to TP/K\to TM\to 0,$$
and a central extension of the structure group
$$1\to Z\to \hat{K}\to K\to 1,$$
one associates a crossed module of topological Lie algebras. For this, 
denoting ${\mathfrak z}$, ${\mathfrak k}$ and $\hat{\mathfrak k}$ the Lie 
algebras of $Z$, $K$ and $\hat{K}$ respectively, we first associate to the 
extension the sequence of vector bundles
\begin{equation}  \label{extension}
0\to M\times{\mathfrak z}\to\widehat{{\rm Ad}}(P)\to{\rm Ad}(P)\to 0,
\end{equation}
where $\widehat{{\rm Ad}}(P)$ is the bundle $(P\times\hat{\mathfrak k})/K$. 
Observe that the adjoint action of $\hat{K}$ on $\hat{\mathfrak k}$ factors 
to an action of $K$ on $\hat{\mathfrak k}$. The bundle $\widehat{{\rm Ad}}(P)$,
which is isomorphic to ${\rm Ad}(\hat{P})$ in case $\hat{P}$ exists, is 
constructed from the ingredients of the problem and thus exists even if 
$\hat{P}$ does not exist. 

We now pass to the spaces of global sections of the above vector bundles. 
For the Atiyah sequence, we get 
${\mathfrak a}{\mathfrak u}{\mathfrak t}(P)=\Gamma(TP/K)={\mathcal V}(P)^K$ 
the Lie algebra of $K$-invariant vector fields on $P$, 
${\mathcal V}(M)=\Gamma(TM)$ the Lie algebra of vector fields on $M$, and 
${\mathfrak g}{\mathfrak a}{\mathfrak u}(P)={\mathfrak n}=\Gamma({\rm Ad}(P))$ 
the {\it gauge Lie algebra}, i.e. the Lie algebra of vertical $K$-invariant 
vector fields. For the sequence (\ref{extension}), we get Lie algebras 
${\mathfrak c}=\Gamma(M\times{\mathfrak z})$, 
${\mathfrak n}=\Gamma({\rm Ad}(P))$ and 
$\hat{\mathfrak n}=\Gamma(\widehat{\rm Ad}(P))$.
All these are given the topology of uniform convergence on compact sets of the 
function and all of its derivatives, and become in this way locally convex 
topological Lie algebras.

One gets a crossed module of topological Lie algebras 
$\mu:\hat{\mathfrak n}\to{\mathfrak a}{\mathfrak u}{\mathfrak t}(P)$ 
by projecting onto ${\mathfrak n}$ and including then ${\mathfrak n}$ into 
${\mathfrak a}{\mathfrak u}{\mathfrak t}(P)$, or, passing to the kernel and 
the cokernel of $\mu$, we get a four term exact sequence
$$0\to{\mathfrak c}\stackrel{i}{\to}\hat{\mathfrak n}\stackrel{\mu}{\to}
{\mathfrak a}{\mathfrak u}{\mathfrak t}(P)\stackrel{\pi}{\to}{\mathcal V}(M)
\to 0.$$
The action $\eta$ of ${\mathfrak a}{\mathfrak u}{\mathfrak t}(P)$ on 
$\hat{\mathfrak n}$ is induced by the derivation action of ${\mathcal V}(P)$
on $\hat{\mathfrak n}\subset{\mathcal C}^{\infty}(P,\hat{\mathfrak k})$. 

The differential $3$-form $\omega_{\rm top\,\,alg}$ in 
$\Omega^3(M,{\mathfrak z})$ is constructed as a cocycle associated to 
this crossed module (cf lemma VI.2 in \cite{Neeb}). Namely, take a 
principal connection $1$-form $\theta\in\Omega^1(P,{\mathfrak k})^K$. 
It serves two purposes: first, it defines a section $\sigma$ of $\pi$, 
second, it gives a projection from 
${\mathfrak a}{\mathfrak u}{\mathfrak t}(P)$ to ${\mathfrak n}$. 
The curvature of $\theta$ is the $2$-form 
$$R_{\theta}=d\theta\,+\,\frac{1}{2}[\theta,\theta]\in\Omega^2(P,
{\mathfrak k})^K.$$
Now, regarding $R_{\theta}$ as an ${\rm Ad}(P)$-valued $2$-form, it can be 
lifted to an $\widehat{\rm Ad}(P)$-valued $2$-form $\Omega$. The section 
$\sigma$, and therefore the connection $\theta$, permit to define an
{\it outer action} 
$$S=\eta\circ\sigma:{\mathcal V}(M)\to{\mathfrak d}{\mathfrak e}
{\mathfrak r}(\hat{\mathfrak n}),$$ 
meaning that the linear map $S$ is a homomorphism when projected to 
${\mathfrak o}{\mathfrak u}{\mathfrak t}(\hat{\mathfrak n})$ and thus gives, 
in particular, an action of ${\mathcal V}(M)$ on ${\mathfrak c}$.
The class $\omega_{\rm top\,\,alg}$ is then just $d_S\Omega$, but interpreted 
as an element of $\Omega^3(M,{\mathfrak z})$. Here $d_S$ means that one takes 
the Chevalley-Eilenberg coboundary with values in $\hat{\mathfrak n}$ as 
if $ S$ were an action:
\begin{eqnarray*}
d_S\Omega(x_1,x_2,x_3)&=&S(x_1)\cdot\Omega(\sigma x_2,\sigma x_3)+
S(x_2)\cdot\Omega(\sigma x_3,\sigma x_1)+S(x_3)\cdot\Omega(\sigma x_1,
\sigma x_2)  \\
&-&\Omega(\sigma [x_1,x_2],\sigma x_3)-
\Omega(\sigma [x_2,x_3],\sigma x_1)-\Omega(\sigma [x_3,x_1],\sigma x_2)
\end{eqnarray*}
for $x_1,x_2,x_3\in{\mathcal V}(M)$. Geometrically speaking, $d_S\Omega$ is 
the covariant derivative $d_S\Omega=d\Omega+S\wedge\Omega$, cf \cite{Neeb} 
section VI. 

\begin{rem}
We consider that we have fixed a decomposition as vector spaces 
$\hat{\mathfrak k}={\mathfrak k}\oplus{\mathfrak z}$, and that with respect 
to this decomposition, the above lift sends $x\in{\mathfrak k}$ to 
$(x,0)\in{\mathfrak k}\oplus{\mathfrak z}$.
Observe that as $d_S\Omega(x_1,x_2,x_3)$ has values in ${\mathfrak z}$, 
$R_{\theta}$ has 
values in ${\mathfrak k}$ and each term 
$S(x_i)\cdot\Omega(\sigma x_j,\sigma x_k)$ has a component in 
${\mathfrak z}$ and one in ${\mathfrak k}$, the sum in the second line 
of the equation for $d_S\Omega(x_1,x_2,x_3)$ together with the components with
values in ${\mathfrak k}$ from the first line must give zero. $d_S\Omega$
is thus just the component with values in ${\mathfrak z}$ of the cyclic sum 
of the terms of the form $S(x_i)\cdot\Omega(\sigma x_j,\sigma x_k)$.
\end{rem}

\section{Crossed modules of Lie algebroids}

We recall in this section the main definitions on Lie algebroids, their 
cohomology and crossed modules, cf \cite{Andr}, \cite{Mack}.

\begin{defi}\label{def:algebroids}
Let $M$ be a fixed manifold. A Lie algebroid $A$ over $M$ is a vector bundle 
$q:A\to M$ together with a morphism of vector bundles $a:A\to TM$ from $M$ to 
the tangent bundle $TM$ of $A$, called the anchor, and an $\R$-linear 
skewsymmetric bracket $[,]:\Gamma A\times\Gamma A\to\Gamma A$ on the space of 
smooth sections $\Gamma A$ of $A$ which satisfies the Jacobi identity, and 
such that

\begin{itemize}
\item[(1)] $[X,fY]\,=\,f[X,Y]+a(X)(f)Y$
\item[(2)] $a([X,Y])\,=\,[a(X),a(Y)]$ (this second point being in fact
a consequence of the first one). 
\end{itemize}
A Lie algebroid $A$ is called transitive if $a$ is fibrewise surjective. In 
this case, $a$ is of constant rank and $TM$ is a quotient bundle of
$A$. 
Given two Lie algebroids $A$ and $A' $  over the same manifold $M$,
a
morphism $\phi:A\to A'$ of vector bundles
such that, first, $q' \circ \phi =q $, second
 $\phi$  induces a Lie algebra homomorphism on the spaces of 
sections, and such that $a'\circ\phi\,=\,a$ is said to be a Lie
algebroid homomorphism over the identity of $M$. 
\end{defi}

We may say simply {\it Lie algebroid morphism},
when it is obvious that it is over the identity of $M$.
 
Note that if the anchor map is trivial at each point of $M$,
a Lie algebroid is precisely what is called a {\it Lie algebra bundle},
i.e. a vector bundle, endowed fiberwise with a Lie algebra structure,
which is required to satisfy the following assumption of regularity: 
for two smooths sections of the vector bundle, taking the bracket
pointwise yields an other smooth section.

The main example of a Lie algebroid and Lie algebroid
morphisms which we have in mind is the Atiyah 
algebroid and its associated sequence of a principal $K$-bundle $P$ on a 
manifold $M$, cf \cite{Mack} p. 97. It is a transitive Lie algebroid.    

\begin{rem}
Regarding a Lie algebra as a Lie algebroid over the point defines a fully 
faithful functor from the category of Lie algebras to the category of Lie 
algebroids. Sending a Lie algebroid to its space of global sections defines 
a functor $\Gamma$ from the category of Lie algebroids to the category of 
(topological) Lie algebras which is an inverse to the previous one when 
restricted to Lie algebroids over the point.
\end{rem} 

\begin{defi}
Let $A$ be a Lie algebroid on $M$, and $E$ be a vector bundle on $M$. A 
representation of $A$ on $E$ is a morphism of Lie algebroids 
$\rho:A\to{\mathcal D}(E)$ where ${\mathcal D}(E)$ is the Lie algebroid of 
first order differential operators $D:\Gamma E\to \Gamma E$ such that there 
is a vector fields $X=a_{\mathcal D}(D)$ with 
$D(fs)\,=\,fD(s)+a_{\mathcal D}(D)(f)s$. 

We shall then say that $E$ is an $A$-module. A morphism of representations 
$\phi:E\to E'$ where $E$, $E'$ are $A$-modules via $\rho$ and $\rho'$, is a 
morphism of vector bundles satisfying $\phi(\rho(X)(s))=\rho'(X)(\phi(s))$ 
for all $s\in\Gamma E$, $X\in\Gamma A$.
\end{defi}

For instance, the trivial bundle $M \times {\mathbb R}^k  \to M$ is an
$A$-module, when the Lie algebroid morphism $\rho$ of the definition
is choosen to be the anchor map (and vector fields are considered as
derivations of the space of $k $-tuples of smooth functions). 

Let us now recall the complex which computes the cohomology of a Lie 
algebroid with values in some representation:

\begin{defi}
Let $A$ be a Lie algebroid and $E$ be an $A$-module, both being vector bundles 
on $M$. The standard complex $C^*(A,E)$ is the complex of vector bundles 
$C^p(A,E)\,=\,\Alt^p(A,E)$ together with the usual differential 
$d:\Gamma C^p(A,E)\to\Gamma C^{p+1}(A,E)$. Cocycles, coboundaries and 
cohomology spaces are defined in the usual way.
\end{defi}

The tangent bundle $TM \to M$ is a Lie algebroid: the anchor map is
the identity and the bracket is the bracket of vector fields. 
With respect to the trivial representation  $M \times {\mathbb R}^k  \to M$,
the space of cochains $C^p(A,E)$ is precisely the space of $k$-tuples
of $p$-forms.
It is a direct verification that the algebroid differential becomes,
under this identification, the de Rham differential. 
In conclusion,

\begin{lem}\label{lem:deRham-algebroid}
The algebroid cohomology of the tangent bundle $TM \to M $ 
taking values in the
trivial module $M \times {\mathbb R}^k  \to M $  is  the de Rham
cohomology tensored by ${\mathbb R}^k $, in short 
$$ H^*(A) = \bigoplus_{i \in  \mathbb
  N} H_{\rm dR}^i(M) \otimes  {\mathbb R}^k .$$
\end{lem}

This lemma stays true for values in a vector bundle with infinite dimensional
fiber, seen as a trivial $TM \to M $-module.

The following definition is taken from \cite{Andr}: 

\begin{defi}
A crossed module of Lie algebroids over the manifold $M$ is a quadruple 
$(L,\mu,A,\rho)$ where $L\to M$ is a Lie algebra bundle, $A\to M$ is a 
transitive Lie algebroid, $\mu:L\to A$ is a morphism of Lie algebroids and 
$\rho:A\to {\mathcal D}(L)$ is a representation of $A$ in $L$ such that 
\begin{itemize}
\item[(a)] $\rho(X)([U,V])\,=\,[\rho(X)(U),V]+[U,\rho(X)(V)]$ for all
  $X\in\Gamma A$ and for all $U,V\in\Gamma L$ and
\item[(b)] $\rho(\mu(V))(W)\,=\,[V,W]$ for all $V,W\in\Gamma L$ and
\item[(c)] $\mu(\rho(X)(V))\,=\,[X,\mu(V)]$ for all $X\in\Gamma A$ and 
$V\in\Gamma L$.
\end{itemize}
\end{defi}

Observe that this notion is actually the notion of a crossed module of 
{\it transitive} Lie algebroids.

According to \cite{Andr}, Section 2, the map $\mu $ is of constant rank.

Let us recall \cite{Mack} p. 272-273 how to associate a $3$-cocycle 
$$\omega_{\rm alg} \in  C^3(\coker(\mu),\ker(\mu))$$ 
to a crossed module of Lie algebroids. Every crossed 
module of Lie algebroids $(L,\mu,A,\rho)$ induces (using $\rho$) a 
{\it coupling}, i.e. a morphism of Lie algebroids 
$\coker(\mu)\,\to\,{\rm Out}{\mathcal D}(L)$. A lift to 
$\nabla:\coker(\mu)\,\to\,{\mathcal D}(L)$ is called a 
{\it Lie derivation law} covering the coupling. $\nabla$ can be taken formally 
to be an action of $\coker\mu$ on $L$, although it is not an action in 
general. The curvature of $\nabla$ takes values in the inner derivations of 
$L$, and can thus be lifted to an alternating vector bundle map 
$\Lambda:(\coker\,\,\mu)\oplus(\coker\,\,\mu)\to L$. The cocycle is now 
obtained as
$$\gamma(X,Y,Z):=\sum_{\rm cycl.}\left(\nabla_X(\Lambda(Y,Z))-\Lambda([X,Y],Z)
\right),$$
i.e. as the formal Chevalley-Eilenberg coboundary operator applied to 
the cochain $\Lambda$ using the ``action'' of $\coker\,\,\mu$ on $L$. 
It is clear that $\nabla$ can be given by a section of the quotient map 
$A\to\coker\,\,\mu$ and the action $\rho$.

Now let us look at the special case which is of interest in this paper: let 
$P$ be a $K$-principal bundle on a manifold $M$, and 
$$1\to Z\to \hat{K}\to K\to 1$$
be a central extension of the structure group. To these data, we attach the 
crossed module of Lie algebroids $\mu:\widehat{{\rm Ad}}(P)\to TP/K$. Here, 
$\widehat{{\rm Ad}}(P)$ denotes the bundle of Lie algebras 
$(P\times\hat{\mathfrak k})/K$, $TP/K$ is the transitive Atiyah Lie algebroid, 
and a representation of $TP/K$ on $\widehat{{\rm Ad}}(P)$ by first order 
differential operators is given by the action of $K$-invariant vector fields 
on $P$ as derivations on functions from $P$ to $\hat{\mathfrak k}$. 
This action makes $\mu$ equivariant and restricts to the action of 
${{\rm Ad}}(P)$ on $\widehat{{\rm Ad}}(P)$ given by the central
extension. 
In this case, $\coker(\mu) $ is the tangent bundle $TM \to M$ 
and $ \ker(\mu)$ is the trivial bundle $M \times {\mathfrak z} \to M
$,
which amounts to the fact that  $\omega_{\rm alg} $ is an element of
$\Omega^3(M,{\mathfrak z}) $, and, by Lemma \ref{lem:deRham-algebroid},
 the class 
$ [\omega_{\rm alg}]$ belongs to $H^3_{\rm dR}(M,{\mathfrak z})=
 H^3_{\rm dR}(M)
 \otimes {\mathfrak z} $.

\begin{theo}
Given a $K$-principal bundle $P$ on a manifold $M$ and a central extension 
$$1\to Z\to \hat{K}\to K\to 1$$
of the structure group, the $3$-cohomology class $[\omega_{\rm top\,\,alg}]$ 
of the crossed module of topological Lie algebras 
$$0\to{\mathfrak c}\stackrel{i}{\to}\hat{\mathfrak n}\stackrel{\mu}{\to}
{\mathfrak a}{\mathfrak u}{\mathfrak t}(P)\stackrel{\pi}{\to}{\mathcal V}(M)
\to 0$$
coincides with the class $[\omega_{\rm alg}]$ associated to the crossed 
module of Lie algebroids 
$$0\to M\times{\mathfrak z}\to\widehat{{\rm Ad}}(P)\to TP/K\to TM\to 0.$$
\end{theo}

\pr 

Choosing a principal connection $1$-form $\theta$ gives a horizontal lifting
for the Lie algebroids, as for the Lie algebras of global sections. $\theta$
gives rise on the one hand to a coupling and a $\nabla$ which corresponds on 
the other hand to the outer action $S=\eta\circ\sigma$ where $\sigma$ is 
also defined by $\theta$. The curvature $R_{\theta}$ of $\theta$ can be 
regarded as an ${\rm Ad}(P)$-valued $2$-form. It is then lifted to an 
$\widehat{\rm Ad}(P)$-valued $2$-form $\Omega$ on $M$ which corresponds to 
$\Lambda$ in the discussion preceeding the theorem. 
The next step is to compute the above mentioned formal Chevalley-Eilenberg 
coboundary of $\Omega$ with values in $\widehat{\rm Ad}(P)$
denoted by $d_S$ in the previous section. Up to identification, we get thus
that for this special crossed module of Lie algebroids, $\omega_{\rm alg}$
has the same expression as $\omega_{\rm top\,\,alg}$.  

\fin

By the above proposition, we have at our disposal all results on Lie 
algebroids in part $2$, ch. $7$ of \cite{Mack}; applying the global section 
functor, they give results on the class of $\omega_{\rm top\,\,alg}$. For 
example, Neeb shows in \cite{Neeb} prop. VI.3 that if an extension 
$\hat{P}$, i.e. a principal $\hat{K}$-bundle $\hat{P}$ such that 
$\hat{P}/Z\cong P$, exists, the class of $\omega_{\rm top\,\,alg}$ is 
trivial. But on the other hand, cor. 7.3.9 on p. 281 in \cite{Mack} shows 
that in case $[\omega_{\rm alg}]=0$, there exists an extension of Lie 
algebroids
$$0\to \widehat{\rm Ad}(P)\to R\to TM\to 0.$$
This is in some sense the converse of Neeb's prop. VI.3. However, we do not 
know how $R$ is related to the existence of a $\hat{K}$-principal 
bundle $\hat{P}$, but we will see that the question of the existence 
of $\hat{P}$ can be solved in terms of Lie groupoids associated to the 
Lie algebroids studied here.

\section{Crossed modules of Lie groupoids}

\begin{defi}
A Lie groupoid consists of two manifolds $\Gamma$ and $M$, together with two 
surjective submersions $s,t:\Gamma\to M$, called the source and the target 
map, and a smooth object inclusion $M\to \Gamma$ such that for composable 
elements (i.e. $g,h\in \Gamma$ such that $t(g)=s(h)$), there is a smooth 
group law having the images of the object inclusion as its unit elements.

A morphism of Lie groupoids ${\cal F}:\Gamma\to \Xi$ where $\Gamma$ is on $M$ 
and $\Xi$ is on $N$, is a pair of smooth maps 
$({\cal F}:\Gamma\to \Xi,f:M\to N)$ intertwining sources and targets such 
that ${\cal F}$ is a homomorphism of groups for composable elements. 
\end{defi}

A Lie groupoid where the source and target maps are equal is what is
called a {\it Lie group bundle}, i.e. a bundle $ \Gamma$ over $M$
where all the fibers admit a
Lie group structure which is smooth in the sense
that taking the pointwise product of two smooth local
sections, or the pointwise inverse of a local smooth section,
yields a smooth local section.

We do not limit ourself to the case of manifolds of finite dimension, but 
$M$ and $H^*_{\rm dR}(M)$ will always be supposed finite dimensional;
as remarked before, Lemma \ref{lem:deRham-algebroid} can be easily adapted 
to (this version of) the infinite dimensional framework. If $\Gamma$ is 
infinite dimensional, one can use the notion of submersion with the help
of an implicit function theorem with parameter like in \cite{Gloe}, as it is
shown explicitly in \cite{NeebWage}.

\begin{defi}
Let $\Gamma$ be a Lie groupoid on $M$, and let $\pi:F\to M$ be a bundle of 
Lie groups on $M$. A representation of $\Gamma$ on $F$ is a smooth map 
$\rho:\Gamma\times_M F\to F$ (where the fiber product is taken with respect 
to the source map $s$), such that 
\begin{itemize}
\item $\pi\circ\rho= t\circ{\rm proj}_1$, with the target map $t$,
\item $\rho$ is an action for composable elements, and
\item $\rho(x):F_{s(x)}\to F_{t(x)}$ is a Lie group isomorphism for all 
$x\in\Gamma$.
\end{itemize}
\end{defi}

\begin{defi}
A crossed module of Lie groupoids ${\cal F}:F\to\Gamma$ is a Lie 
groupoid $\Gamma$ on $M$, a Lie group bundle $F$ on $M$, a morphism of Lie 
groupoids ${\cal F}:F\to\Gamma$ (with $f$ the identity), and a 
representation $\rho$ of $\Gamma$ on $F$ such that 
\begin{itemize}
\item ${\cal F}(\rho(x,y))=x{\cal F}(y)x^{-1}$ for all 
$(x,y)\in\Gamma\times_MF$, 
\item $\rho({\cal F}(x),y)=xyx^{-1}$ for all $x,y\in F$ avec $\pi(x)=\pi(y)$, 
and
\item the image of ${\cal F}$ is a closed regularly embedded submanifold of 
$\Gamma$.
\end{itemize}
\end{defi}

We will only work here with crossed modules of Lie groupoids 
${\cal F}:F\to\Gamma$ such that $\coker\,\,{\cal F}=M\times M$  (the pair 
groupoid)  and 
$\ker({\cal F}) \simeq M \times Z$ for some Abelian Lie group $Z$.
As a crossed module of transitive Lie algebroids is in a sense a
crossed module of $TM$ by a trivial ${\mathfrak z}$-bundle, the crossed modules
of Lie groupoids we discuss here are crossed modules of $M\times M$ by a
trivial Lie group bundle $Z\times M$. 
In other words, we restrict ourself to the case of crossed modules over  
{\it transitive} Lie groupoids.

An important point is that, in our case, $\Gamma $ and $F$ may be infinite
dimensional, while $\coker \, \, {\cal F}$ is finite dimensional.

Suppose given a $K$-principal bundle $P$ and a central extension $\hat{K}$ of
the structure group $K$. Note that $K$ acts by conjugation, not only on 
$K$ itself, but also on $\hat{K}$.
Denote by $P_K(\hat{K})={\mathcal C}^{\infty}(P,\hat{K})^K$ the space of 
$K$-equivariant smooth maps from $P_m$ to $\hat{K}$, where $P_m$ is the 
fiber over an arbitrary $m \in M$.
$P_K(\hat{K})$ is naturally a group bundle over $M$, and there 
is a natural groupoid homomorphism
${\cal F}: P_K(\hat{K}) \to (P \times P)/K $ mapping
$\phi\in{\mathcal C}^{\infty}(P_m,\hat{K})^K $ to 
$\overline{(p\cdot\nu(\phi(p)),p)}$
where  $p$ is an arbitrary element of $P_m$, $\nu$ is the map from 
$\hat{K}$ to $P$, and the bar means the class in $(P\times P)\,/\,K$.
It is easy to check that the groupoid homomorphism
${\cal F}$ is a crossed module. Conversely, we have:

\begin{prop}  \label{prop1}
Let ${\cal F}:F\to\Gamma$ be a crossed modules of Lie groupoids 
such that $\coker\,\,{\cal F}=M\times M$ and such 
that the kernel of ${\cal F}$ is trivial: $\ker\,{\cal F}=Z\times M$.

Fix a point $m \in M$, and define the Lie group $K$ to be the
quotient of $\hat{K}:=F_m$ by $Z$. Then there exists a $K$-principal bundle 
$P$ such that ${\cal F}:F\to\Gamma$ is isomorphic to the crossed module
$$0\to M\times Z\to P_K(\hat{K})\to (P\times P)/K\to M\times M\to 0.$$
\end{prop}

\pr
We fix a point $m \in M$, and define the Lie group $K$ to be the
quotient of $\hat{K}:=F_m$ by $Z$. 
Denote by $s$ and $t$ the source and target maps of $\Gamma$. 
The space $t^{-1}(m)$ is a submanifold of $\Gamma$
(because the target map is a submersion) on which $K$
acts freely by right multiplication.
Now, the source map $s:t^{-1}(m)\to M$ is a surjective submersion 
onto $M$ as $\coker({\cal F})=M\times M$, and the fibers are precisely given
the right $K$-action. Hence $ t^{-1 } (m) \to M$ is  a principal
$K$-bundle that we denote by $P$. 

It is easy to check that the groupoid $\Gamma$ is isomorphic to the
Atiyah groupoid $(P \times P)/K $. The isomorphism $\Psi$ is as follows. 
Any element of $\Gamma $ can be written in the form $\gamma = \gamma_1 \cdot
\gamma_2 $
with $\gamma_1,\gamma_2^{-1} \in t^{-1}(m) $. Define 
$\Psi(\gamma)  $ to be $  \overline{(\gamma_1,\gamma_2^{-1})} \in  
(P \times P)/K$. 

Let us also define a map 
$\Phi:F\to P_K(\hat{K})={\mathcal C}^{\infty}(P,\hat{K})^K$. To an $f\in F$, 
we associate the map $(x\mapsto \rho(x,f))$, where $\rho$ is the action of
$\Gamma$ on $F$, given in the data of the crossed module, and let us recall
that the map ${\mathcal G}:P_K(\hat{K})\to (P\times P)\,/\,K$ of the Atiyah
crossed module corresponding to the central extension of the structural 
group $K$ of $P$ to $\hat{K}$, is given in our context by
$$\phi\mapsto {\mathcal G}(\phi):=\overline{(p\cdot_{\Gamma}{\mathcal F}
(\phi(p)),p)},$$
where $p\in P$ is an arbitrary point, $\cdot_{\Gamma}$ is the 
multiplication in $\Gamma$. 

Let us show that the square 

\vspace{.5cm} 
\hspace{3cm}
\xymatrix{
P_K(\hat{K})  \ar[r]_{\mathcal G}  &   (P\times P)\,/\,K \\
F \ar[u]_{\Phi} \ar[r]_{\mathcal F} & \Gamma \ar[u]_{\Psi}}
\vspace{.5cm}
 
\noindent is commutative. Indeed, by the axioms of a crossed module, we have a 
commutative square
 
\vspace{.5cm} 
\hspace{3cm}
\xymatrix{
F_n  \ar[r]_{{\mathcal F}_n}  \ar[d]_{\rho(p,1)} &  \Gamma_n^n \\
F_m \ar[r]_{{\mathcal F}_m} & \Gamma_m^m \ar[u]_{{\rm Conj}(p^{-1})}}
\vspace{.5cm}

\noindent where $n=t'(p)$ and $m=s'(p)$, $t'$ and $s'$ being the source and 
target maps 
of $P$ (i.e. the restrictions of those of $\Gamma$ to $P\subset\Gamma$). 
In $n$, we have ${\mathcal F}(f)=\gamma_1\gamma_2^{-1}$ with 
$\gamma_1,\gamma_2\in\Gamma_n^m$ and $f\in F_n$. Now compute
$$\Psi({\mathcal F}(f))\,=\,\overline{(\gamma_1,\gamma_2)}\,=\,
\overline{({\mathcal F}(f)p,p)}\,=\,\overline{(pp^{-1}{\mathcal F}(f)p,p)}
\,=\,\overline{(p{\mathcal F}(\rho(p,f)),p)}.$$

\fin

There is a standard way to associate to a crossed module of Lie groupoids
a characteristic class $\omega_{\rm grp}$ (cf \cite{Mack1} p. 197 or 
\cite{Andr} p. 13): 
let us choose a covering 
(trivializing the principal bundle $P$ described above) ${\cal U}=\{U_i\}$ 
on the manifold 
$M$, and denote as usual $U_{ij}=U_i\cap U_j$, $U_{ijk}=U_i\cap U_j\cap U_k$. 
The principal bundle $P$ is given by transition functions 
$g_{ij}:U_{ij}\to K$, which we may as well regard as a \v{C}ech $1$-cocycle. 
Lift the functions $g_{ij}$ to functions $\hat{g}_{ij}$ with values in 
$\hat{K}$, and denote by $\hat{g}_{ij}\hat{g}_{jk}\hat{g}_{ki}=h_{ijk}$ 
their default for forming a \v{C}ech $1$-cocycle. 
$h_{ijk}$ has values in the sheaf $\underline{Z}$ of differentiable 
$Z$-valued functions, as the $g_{ij}$ do form a cocycle. 
The \v{C}ech $2$-cocycle $h_{ijk}$ with values in $\underline{Z}$
is by definition $\omega_{\rm grp}$. Its class will also be regarded as 
a class in $H^3(M,\pi_1(Z))$, provided that $Z$ is a connected regular 
abelian Lie group, 
according to the isomorphism 
$$H^2(M,\underline{Z})\,\cong\,H^3(M,\pi_1(Z)),$$
which stems from the fact that for a connected regular abelian Lie group
$$1\to \pi_1(Z)\to{\mathfrak z}\to Z\to 1$$ 
is an exact sequence of groups, where the abelian group ${\mathfrak z}$ is
the Lie algebra of $Z$, and that the sheaf $\underline{\mathfrak z}$ is fine.

The last step is to consider the image of $H^3(M,\pi_1(Z))$ in the de
Rham cohomology group $H^3_{\rm dR}(M)\otimes {\mathfrak z}$ obtained by 
composing the map  $H^3(M,\pi_1(Z)) \to H^3(M,{\mathbb R}) \otimes 
{\mathfrak z}$ coming from  the inclusion $\pi_1(Z) $ into
${\mathfrak z} $ with the \v{C}ech-de Rham isomorphism 
$H^3(M,{\mathbb R}) \, \simeq \, H^3_{\rm dR}(M) $.

The following theorem is the main result of our paper:

\begin{theo}  \label{mainthm}
Given a $K$-principal bundle $P$ on a manifold $M$ and a central extension 
$$1\to Z\to \hat{K}\to K\to 1$$
of the structure group, the $3$-cohomology class $[\omega_{\rm top\,\,alg}]$ 
of the crossed module of topological Lie algebras 
$$0\to{\mathfrak c}\stackrel{i}{\to}\hat{\mathfrak n}\stackrel{\mu}{\to}
{\mathfrak a}{\mathfrak u}{\mathfrak t}(P)\stackrel{\pi}{\to}{\mathcal V}(M)
\to 0$$
defines the same de Rham cohomology class as the $3$-cocycle 
$\omega_{\rm grp}$ associated to the crossed module of Lie groupoids 
$$0\to M\times Z\to P_K(\hat{K})\to (P\times P)/K\to M\times M\to 0.$$
\end{theo}

\pr We use the previously established notations. Let us choose a
connection on the principal bundle $P \to M $ given, in local
trivializing coordinates $\{U_i\}_{i \in I}$,
by  a family $\theta_i$ of ${\mathfrak k} $-valued $1$-forms on $U_i$. 
As usual, the relation  
$\theta_j={\rm Ad}_{g_{ji}}\theta_i+g_{ji}^{-1}dg_{ji}$ expresses how to pass
from $\theta_i$ to $\theta_j$ by gauge transformation. Lifting to 
$\hat{\mathfrak k}$, we get 
$\hat{\theta}_j-{\rm Ad}_{\hat{g}_{ji}}\hat{\theta}_i=\hat{g}_{ij}^{-1}
d\hat{g}_{ij}-\alpha_{ij}$ with $\alpha_{ij}\in{\mathfrak z}$
(remember that we fixed a vector space decomposition
$\hat{\mathfrak k}={\mathfrak k}\oplus{\mathfrak z}$). 
We get furthermore for the curvature as usual
$R_i=d\theta_i+\frac{1}{2}[\theta_i,\theta_i]$,
$R_i={\rm Ad}_{g_{ij}}R_j$ and  
$dR_i=[\theta_i,R_i]$, and for the lifted curvature
$\hat{R}_i=d\hat{\theta}_i+\frac{1}{2}[\hat{\theta}_i,\hat{\theta}_i]$,
$\hat{R}_i={\rm Ad}_{\hat{g}_{ij}}\hat{R}_j+d\alpha_{ij}$.

The \v{C}ech $2$-cocycle $h_{ijk}$ represents by definition the class 
$\omega_{\rm grp}$. In order to establish the link with the class 
$\omega_{\rm top\,\,alg}$, we translate \v{C}ech cocycles into differential 
forms via the \v{C}ech-de Rham bicomplex: a straight forward computation,
using that $h_{ijk}$ and $\alpha_{ij}$ have values in $Z$, gives

\begin{equation}   \label{****}
\alpha_{ij}-\alpha_{kj}+\alpha_{ki}\,=\,h_{ijk}^{-1}dh_{ijk}.
\end{equation}

Now choose a form
$\omega_i\in\Omega^2(U_i,\hat{\mathfrak k})$ which has the right 
transformation 
behaviour $\omega_j={\rm Ad}_{\hat{g}_{ij}}\omega_i$ and which has 
$R_i$ as its 
component in ${\mathfrak k}$ and a form $\omega^{\mathfrak z}$ as its 
component in ${\mathfrak z}$. Observe that the ${\mathfrak z}$-component
$\omega^{\mathfrak z}$ is globally defined. It is then the form
$\beta_i:=\hat{R}_i-\omega_i\in\Omega^2(U_i,{\mathfrak z})$ which satisfies
$\beta_i-\beta_j=d\alpha_{ij}$.

The relation (\ref{****}) means that $(h,\alpha)$ forms a \v{C}ech $2$-cocycle
with values in the complex of sheaves
$$d\log:\Omega^0_M(Z)\to\Omega^1_M({\mathfrak z}).$$

Then $d\beta_i$, which is just the component in ${\mathfrak z}$ of 
$d\hat{R}_i=[\hat{\theta}_i,\hat{R}_i]$ up to an exact form,  
is exactly the expression of a 
representative of $\omega_{\rm top\,\,alg}$ (cf remark $1$). 

The general procedure means here that $(h,\alpha,\beta)$ forms a \v{C}ech 
$2$-cocycle with values in the complex of sheaves

\begin{equation}  \label{*****}
\Omega^0_M(Z)\stackrel{d\log}{\to}\Omega^1_M({\mathfrak z})\stackrel{d}{\to}
\Omega^1_M({\mathfrak z}).
\end{equation}

The rest of the proof of the theorem is inspired by the proof of proposition 
$4.2.7$ in \cite{Bryl} p. $174$:
By corollary $2$ in the appendix, $(h,\alpha,\beta)$
defines a cohomology class in the smooth Deligne cohomology group
$H^3(M,\pi_1(Z)(3)^{\infty}_D)$. Furthermore, the class is sent to
the \v{C}ech-connecting homomorphism of $g_{ij}$, i.e. $h_{ijk}$, under the 
map in cohomology from $H^3(M,\pi_1(Z)(3)^{\infty}_D)$ to 
$H^3(M,\pi_1(Z))\cong H^2(M,\underline{Z})$, induced by the canonical 
projection from the complex of sheaves $\pi_1(Z)(3)^{\infty}_D$ to its first 
term (which is the constant sheaf $\pi_1(Z)$), see proposition $2$ in the
appendix. On the other hand, the Deligne cohomology class maps to 
$[d\beta]$, and thus to $[\omega_{\rm top\,\,alg}]$, under the map of 
complexes of sheaves given by prolongation of the complex one step further,
see proposition $3$ of the appendix. It is therefore clear that the image of 
$[h_{ijk}]$ or $[\omega_{\rm grp}]$ under the map 
$$H^2(M,\underline{Z})\cong H^3(M,\pi_1(Z))\to H^3(M,{\mathfrak z})$$
is the class $[d\beta]$ or $[\omega_{\rm top\,\,alg}]$.\fin
 
\begin{rem}
In conclusion, the Lie groupoid class $\omega_{\rm grp}$ and the Lie 
algebroid class $\omega_{\rm alg}$ coincide up 
to torsion (i.e. have the same image in $H^3(M,{\mathfrak z})$)
in the above context. 
\end{rem}

\begin{cor}
Let ${\mathcal F}:F\to\Gamma$ be a crossed module of Lie groupoids such that
$\coker\,{\mathcal F} = M\times M$ and such that the Lie group bundle
$\ker\,{\mathcal F}=:Z\times M$ is trivial. Set as in proposition \ref{prop1}
$K:=F_x/Z$, and $P$ for the principal $K$-bundle $t^{-1}(m)\to M$.

Then the $3$-cohomology class $[\omega_{\rm top\,\,alg}]$ 
of the crossed module of topological Lie algebras 
$$0\to{\mathfrak c}\stackrel{i}{\to}\hat{\mathfrak n}\stackrel{\mu}{\to}
{\mathfrak a}{\mathfrak u}{\mathfrak t}(P)\stackrel{\pi}{\to}{\mathcal V}(M)
\to 0$$
defines the same de Rham cohomology class as the $3$-cocycle 
$\omega_{\rm grp}$ associated to the crossed module of Lie groupoids 
$$0\to M\times Z\to P_K(\hat{K})\to (P\times P)/K\to M\times M\to 0.$$
\end{cor} 

\pr This follows immediately from theorem \ref{mainthm} and proposition 
\ref{prop1}. \fin


\section{Finite dimensional structure group}

In this section, we look at the special case of a finite dimensional 
structure group $K$. Let us start with $K=PU_n$, the projective unitary group 
of an $n$-dimensional complex vector space. $K$-principal bundles on $M$
which are non-trivial in the sense that they cannot be lifted to a principal
$U_n$-bundle define elements of the {\it Brauer group} $Br(M)$, cf
\cite{Grot2}. A theorem of Serre determines $Br(M)$:

\begin{theo}[Serre] 
On the manifold $M$, $Br(M)$ can be identified with the torsion subgroup of 
the sheaf cohomology group $H^2(M,\underline{\C^*})\cong H^3(M,\Z)$.
\end{theo}

The identification is given by the obstruction class $[\omega_{\rm grp}]$
which measures the obstruction for a given $PU_n$-principal bundle $P$ to be
lifted to a $U_n$-principal bundle. By Serre's theorem together with 
theorem $2$, we arrive thus at the following conclusion: 

\begin{cor}
Neeb's class $[\omega_{\rm top\,\,alg}]$ associated to the problem of lifting 
a given $PU_n$-principal bundle $P$ to a $U_n$-bundle is always zero.
\end{cor}

We can go a step further. Murray shows in \cite{Mur} \S 13 that the class 
$[d\beta_i]=[\omega_{\rm top\,\,alg}]$ (in the notation of the proof of 
theorem $2$) is trivial in case the structure group is finite dimensional and 
structure group and base space are $1$-connected (this condition is not 
written in \S 13 of \cite{Mur}, but as he refers to \cite{GLSW}, the condition 
seems inevitable).  

\begin{cor}
Let $K$ be a finite dimensional $1$-connected group, $P$ a $K$-principal 
bundle with $1$-connected base, and $\hat{K}$ a central extension of $K$ by
$Z$.

Then Neeb's class $[\omega_{\rm top\,\,alg}]$ associated to the problem of 
lifting $P$ to a $\hat{K}$-bundle $\hat{P}$ such that $\hat{P}/Z\cong P$ is 
zero.
\end{cor}  

We believe that the study of the triviality of Neeb's class in the general case
is an interesting subject.

On the other hand, in some cases, the significance of the (torsion) class 
$[\omega_{\rm grp}]$ is well known; see examples in \cite{Mack1} pp. 206-207. 
To cite just one
example, for $K=SO(n)$ and a principal $K$-bundle $P$ on $M$, $n\geq 3$, 
the obstruction class in $H^2(M,\Z_2)$ for the existence of a $Spin(n)$-bundle
lifting $P$ is the second Stiefel-Whitney class and gives a $2$-torsion element
in $H^3(M,\Z)$ under the Bockstein map.    

\section{Existence of principal bundles}

We will now use the previously established results to answer the question
concerning the existence of a principal $\hat{K}$-bundle $\hat{P}$ such that
$\hat{P}\,/\,Z\cong P$. We use for this the general theory of fiber spaces
set up by Grothendieck \cite{Grot1}, and we will get back the classification 
of principal bundles via cohomology classes of Lie groupoids due to 
Mackenzie \cite{Mack1}. 

Indeed, the short exact sequence of groups
$$1\to Z\to \hat{K}\to K\to 1$$
induces (regardless of the finiteness of the dimension of $K$) an exact 
sequence
$$H^1(M,\underline{\hat{K}})\to H^1(M,\underline{K})\to H^2(M,\underline{Z}).$$
The obstruction class $[\omega_{\rm grp}]$ is nothing else than the image of 
the class of a given $K$-principal bundle $P$ under the connecting 
homomorphism, i.e. the map on the right, and the exactness of the sequence 
means that there exists a $\hat{K}$-principal bundle $\hat{P}$ extending $P$ 
if and only if $[\omega_{\rm grp}]=0$: 

\vspace{.5cm} 
\hspace{2cm}
\xymatrix{
H^1(M,\underline{\hat{K}})  \ar[r]  & H^1(M,\underline{K}) \ar[r] &   
H^2(M,\underline{Z}) \\
[\hat{g}_{ij}]  \ar@{|-->}[r] & [g_{ij}] \ar@{|->}[r] & 
[h_{ijk}=\omega_{\rm grp}]}
\vspace{.5cm} 

This reasoning gives back theorem $3.4^{\scriptstyle \prime}$ in \cite{Mack1}
which shows the existence of a principal bundle $\hat{P}$ if and only if 
$[\omega_{\rm grp}]=0$. By Grothendieck's theory, we get of course a 
\v{C}ech cocycle $\hat{g}_{ij}$ of transition functions which is only 
continuous, but by \cite{Wock}, we may choose a smooth representative of 
the same equivalence class.

Up to torsion, we get thus an equivalence between the existence of a principal
$\hat{K}$-bundle $\hat{P}$ such that $\hat{P}/Z\cong P$ and the vanishing of
Neeb's class $[\omega_{\rm top\,\,alg}]$.

\begin{rem}
Actually, in the context of lifting the structure group to a central extension,
we prescribe not only the outer
action of the Lie groupoid (or principal bundle) whose existence we study, 
but also the Lie algebroid (i.e. $P_K(\hat{\mathfrak k})$). This is the reason
why the \v{C}ech cohomology class $\omega_{\rm grp}$ takes values in a sheaf of
germs of constant maps, cf \cite{Mack1} p. $203$.
\end{rem}

\section{Gerbes}

In this section, we reinterprete our results in terms of gerbes and note
that crossed modules of Lie groupoids and gerbes are related to similar
cohomological problems. We then explore a more direct link between these
two objects.

The theory of gerbes and stacks is rather heavy from the point of 
view of definitions, so we content ourselves with intuitive explanations for
the objects of study in this section, in order not to blow up this paper.

{\it Differentiable stacks} are more general objects than differentiable 
manifolds; their notion grew out of Grothendieck's attempt to characterize the
functor of points of an algebraic variety by abstract conditions. A stack is
in particular a sheaf of categories on a site, but here, our stacks will always
be sheaves on the manifold $M$. The guiding example for us here is the stack 
given by the local liftings of the structure group of a $K$-principal bundle
$P$ from $K$ to a $Z$ central extension $\hat{K}$ of $K$, cf \cite{Bryl} 
pp. 171-172. On an open set $U$ of $M$, the associated category $C_U$ consists 
of $\hat{K}$-principal bundles $\hat{P}$ on $U$ together with an isomorphism
of $K$-principal bundles $f:\hat{P}/Z\cong P$. Morphisms in $C_U$ are bundle
morphisms which commute with the isomorphisms $f$. Note that up to refinement 
of $U$, the category $C_U$ is never empty (for non-empty $U$), and any two 
objects of $C_U$ are isomorphic. These are the basic properties of a
sheaf of groupoids which is a {\it gerbe}.

The functor of points associates to each Lie groupoid a (pre)stack. Up to 
Morita equivalence, this gives a fully faithful functor from Lie groupoids to
differentiable stacks. In this sense a Morita equivalence class of Lie 
groupoids defines a differentiable stack. 

The manifold $M$, seen as a Lie groupoid, is Morita equivalent to the groupoid
associated to an open covering $\coprod_{i\in I}U_i$ of $M$. Now given a class
$\omega\in H^3(M,\Z)$, there is a central extension of Lie groupoids 
(see \cite{TuXu} p. $863$)

\vspace{.5cm} 
\hspace{3cm}
\xymatrix{
R_{\omega}  \ar[r]  & 
\coprod_{i,j\in I}U_{ij} \dar[r] &
\coprod_{i\in I}U_i}
\vspace{.5cm}

On the other hand, the central extension $1\to \C^*\to U({\cal H})\to 
PU({\cal H})\to 1$ of the projective unitary group $PU({\cal H})$ of an
infinite dimensional
separable Hilbert space ${\cal H}$ gives rise to an extension of Lie groupoids

\vspace{.5cm} 
\hspace{3cm}
\xymatrix{
U({\cal H})  \ar[r]  & 
PU({\cal H}) \dar[r] &
{\rm pt}}
\vspace{.5cm}

\noindent over the point ${\rm pt}$. Now it is shown in {\it loc. cit.} 
that there exists
a morphism of stacks from the first extension of Lie groupoids to the second, 
and that this implies the existence of a $PU({\cal H})$-principal bundle $P$ 
on $M$ whose obstruction class is $\omega$.      

All this can be summarized by the morphism $\Phi$
$$\Phi:H^1(M,\underline{K})\times H^2(K,\underline{S^1})\to 
H^2(M,\underline{S^1})$$
to be found in \cite{TuXu} p. 860. Choosing the particular $S^1$-central 
extension of $PU{\cal H}$ given by $U{\cal H}$ (for an infinite dimensional 
separable Hilbert space ${\cal H}$) gives an isomorphism
$$\Phi':H^1(M,{\underline PU{\cal H}})\to H^2(M,\underline{S^1}),$$
see \cite{TuXu} p. 861. We see thus that $PU{\cal H}$-principal bundles (LHS) 
correspond in a one-to-one manner to integral $3$-cohomology of $M$ (RHS) 
via the obstruction mechanism. 
Note that the proofs in this section of \cite{TuXu} 
do certainly not extend to an infinite dimensional center $Z$, and that it 
is not even clear how they extend to a non-compact group $Z$. 

In the framework of \cite{TuXu}, the space $H^2(M,\underline{S^1})$ is 
interpreted as the space of equivalence classes of gerbes; the natural gerbe
associated to our initial problem of lifting the structure group of a principal
bundle to a central extension is the one described in the beginning. The 
above mentioned central extension of Lie groupoids also gives rise to a gerbe.
The {\it band} of a gerbe is the sheaf of groups given by the outer action
(see \cite{LSX} p.15 def. 3.5). In the above central extension, the band 
happens to be trivial.

On the
other hand, we started our study from equivalence classes of crossed modules 
of Lie groupoids. We are therefore led to search a direct link, i.e. one which 
does not pass by cohomology classes, between gerbes and crossed modules of 
Lie groupoids.

\begin{theo}
There is a one-to-one correspondence between gerbes with abelian, 
trivialized band on a connected manifold $M$, and crossed modules of Lie 
groupoids with trivial kernel and cokernel $M\times M$, which induces 
the above one-to-one correspondence between cohomology classes when passing to 
equivalence classes.
\end{theo}

\pr
We will first describe a map $\phi$ associating to a crossed module 
$(F,\delta,\Gamma)$ of Lie groupoids 
$$1\to M\times Z\to F\to \Gamma\to M\times M\to 1$$
a gerbe. $\phi$ will respect cohomology classes in 
$H^2(M,\underline{Z})$.
Recall that $M$ is connected. Call $H$ the fiber $F_{x_0}$ of the Lie group 
bundle $F$ over $x_0\in M$. As described in \cite{Andr} p. 13, one may 
associate to $(F,\delta,\Gamma)$ a cocycle of transition functions 
$s_{ij}:U_{ij}\to\delta(H)$. Denote by  
$\hat{s}_{ij}:U_{ij}\to H$ a lift of $s_{ij}$ with values in $H$. 
We now define a gerbe as a 
sheaf of categories on $M$ by taking for $C_U$ (for an open set $U$ of $M$) 
the category of $H$-principal bundles, given for example by a \v{C}ech cocycle
$\hat{s}_{ij}:U_{ij}\to H$, together with an isomorphism of 
$\delta(H)$-principal bundles $(\delta\circ\hat{s}_{ij})\cong (s_{ij})$. By
the general theory, it is clear that this defines a gerbe 
$C_{(F,\delta,\Gamma)}$ on $M$. We define thus $\phi$ by
$$\phi(F,\delta,\Gamma)\,=\,C_{(F,\delta,\Gamma)}.$$
It is also clear from \cite{Bryl} p. 172 that $\phi$ respects cohomology 
classes.     

Let us now define a map $\psi$ in the reverse direction: a gerbe 
${\mathfrak G}$ (with trivialized abelian band) on $M$ comes together with an 
identification of the sheaf of (abelian) groups of automorphisms of the 
objects of the local categories $C_U$; let this sheaf be $\underline{Z}$ 
and the abelian group be $Z$.
Choosing a local section $s:U\to{\mathfrak G}$ (whose existence is due to
the local existence of objects in $C_U$ - this is one of the 
axioms of a gerbe), 
this identification implies that $U\times_{\mathfrak G}U\to U$ is a locally 
trivial principal bundle $\hat{P}_U$, cf \cite{Hein} rem. 5.2.
The local isomorphy of any 
two objects in $C_U$ (which is the other axiom of a gerbe) implies that the 
fiber of $U\times_{\mathfrak G}U\to U$ is a group. Denote it by $\hat{K}$.   
Now the locally defined $\hat{P}_U/Z=:P_U$ form a globally defined 
$K$-principal bundle $P$ on $M$, because the ${\rm Aut}(C_U)=Z$-valued 
\v{C}ech cochain defined by the local sections becomes a cocycle for $P$. 
This defines as in the beginning of section $4$ a crossed module of 
Lie groupoids $(F,\delta,\Gamma)$ of $M\times M$ by $\underline{Z}$. We set
$$\psi({\mathfrak G})\,=\,(F,\delta,\Gamma).$$

$\phi$ and $\psi$ are mutually inverse bijections, descending to bijections
of equivalence classes. In this sense, gerbes and crossed modules of Lie 
groupoids are the same objects.\fin

\begin{rem}
Note that this construction does not involve the construction of a principal 
bundle $P$ for a given integral cohomology class in \cite{TuXu}. Indeed, 
passing by their construction, one could show a relation between gerbes and
crossed modules only for $K=S^1$.
\end{rem}


\begin{appendix}

\section{Smooth Deligne cohomology}

In this appendix, we present the definitions and results from smooth 
Deligne cohomology which we need in section $4$. Our main reference is 
\cite{Bryl} Ch. 1.5, while our definition differs slightly from his.

Deligne cohomology in our sense is the hypercohomology of truncated 
complexes of sheaves. The complex we consider in section $4$, denoted by 
$\pi_1(Z)(3)^{\infty}_D$, is
$$\pi_1(Z)\stackrel{\rm incl}{\to}\Omega^0_M({\mathfrak z})
\stackrel{d}{\to}\Omega^1_M({\mathfrak z})\stackrel{d}{\to}
\Omega^2_M({\mathfrak z}).$$
The notation $\pi_1(Z)(3)^{\infty}_D$ indicates that the $0$th term is 
the simple (or constant) sheaf $\pi_1(Z)$, and that the complex is 
truncated at the $3$rd place. It will be related to the complex
$$\Omega^0_M(Z)\stackrel{d\log}{\to}\Omega^1_M({\mathfrak z})\stackrel{d}{\to}
\Omega^2_M({\mathfrak z}).$$
Let us explain the ingredients: $Z$ is here some (possibly infinite 
dimensional) Lie group, ${\mathfrak z}$ is 
its Lie algebra. The only property from the setting of infinite dimensional 
Lie groups we might choose is that the composition of smooth maps must be 
smooth. The main geometric property of $Z$ and ${\mathfrak z}$ we use is that
${\mathfrak z}$ is a universal covering space for $Z$, in particular
$$1\to\pi_1(Z)\to{\mathfrak z}\stackrel{\pi}{\to} Z\to 1$$
is an exact sequence of groups. In the above complex, $\Omega^0_M(Z)$ denotes 
the sheaf of smooth maps on $M$ with values in $Z$, and 
$\Omega^i_M({\mathfrak z})$ for $i=1,2$ have a similar meaning. $d\log$ is the
{\it logarithmic derivative}, i.e. $d\log f=f^{-1}df$ for $f\in \Omega^0_M(Z)$.
$d\log f$ for $x\in M$ is thus the composition 
$$L_{f^{-1}(x)}\circ d_xf:T_xM\to T_{f(x)}Z\to T_eZ={\mathfrak z}.$$

\begin{prop}
The following square is commutative:

\hspace{3cm}
\xymatrix{
\Omega^0_M({\mathfrak z}) \ar[r]^{d} \ar[d]^{\pi^*} & 
\Omega^1_M({\mathfrak z}) \ar[d]^{\id} \\
\Omega^0_M(Z) \ar[r]^{d\log} & \Omega^1_M({\mathfrak z}) 
}\vspace{.5cm}
\end{prop} 

\pr It is enough to have $d(\pi\circ f)=(\pi\circ f)\cdot df$ for each 
$f\in \Omega^0_M({\mathfrak z})$ (where $\cdot$ means once again the left 
translation of $df$ in $Z$). 

Observe that for $a\in{\mathfrak z}$, the composition
$$L_{\pi(a)^{-1}}\circ d\pi(a):{\mathfrak z}\to T_{\pi(a)}Z\to
T_eZ={\mathfrak z}$$
is the identity, as $\pi:{\mathfrak z}\to Z$ is the universal covering.
This gives $d\pi(a)=L_{\pi(a)}$ for all $a\in{\mathfrak z}$.

Therefore, $d(\pi\circ f)=(d\pi\circ f)\cdot df=(\pi\circ f)\cdot df$
as claimed.\fin

\begin{cor}
We have a quasi-isomorphism of complexes of sheaves:

\hspace{1cm}
\xymatrix{
\pi_1(Z)  \ar[r]^{{\rm incl}} \ar[d]^{0}  &
\Omega^0_M({\mathfrak z}) \ar[r]^{d} \ar[d]^{\pi^*} & 
\Omega^1_M({\mathfrak z}) \ar[r]^{d} \ar[d]^{\id}  &
\Omega^2_M({\mathfrak z}) \ar[d]^{\id} \\
0 \ar[r]^{0} & 
\Omega^0_M(Z) \ar[r]^{d\log} & 
\Omega^1_M({\mathfrak z}) \ar[r]^{d} &
\Omega^2_M({\mathfrak z})  }\vspace{.5cm}
\end{cor}

\begin{prop}
Let $M$ be a smooth paracompact manifold of dimension $n$, and let
$$0\to{\cal F}^0\to{\cal F}^1\to\ldots\to{\cal F}^n\to 0$$ 
be a complex of sheaves. Denote by 
$H^q(M,{\cal F}^0(p)^{\infty}_D)$ the $q$-th Deligne cohomology group, i.e. the
$q$th hypercohomology group of the truncated sheaf complex 
$\{{\cal F}^i\}_{0\leq i\leq p-1}$. Then the projection of 
$\{{\cal F}^i\}_{0\leq i\leq p-1}$ onto ${\cal F}^0$ induces an epimorphism
$$H^q(M,{\cal F}^0(p)^{\infty}_D)\,\to\,H^q(M,\underline{\cal F}^0)$$
for each $q\geq p$.
\end{prop}

\pr This proposition is analoguous to a part of theorem 1.5.3 in \cite{Bryl}
p. 48 and is shown in the same way.\fin
     
\begin{prop}
In the same setting as the previous proposition, the one-step-extension of
the truncated complex 
$\{{\cal F}^i\}_{0\leq i\leq p-1}$ 

\hspace{1cm}
\xymatrix{
{\cal F}^0  \ar[r] \ar[d]^{0}  &
{\cal F}^1 \ar[r] \ar[d]^0 & 
\ldots \ar[r] \ar[d]^0  &
{\cal F}^{p-1} \ar[d] \\
0 \ar[r]^0 & 
0 \ar[r]^0 & 
\ldots \ar[r] &
{\cal F}^p  }\vspace{.5cm}

induces a map
$$H^p(M,{\cal F}^0(p)^{\infty}_D)\,\to\,{\cal F}^p(M).$$
\end{prop}

\end{appendix}

\end{document}